\documentclass[12pt,a4,amsmath, amsthm, amssymb]{amsart}
\usepackage{amsmath,amssymb,amsthm}

\newtheorem{Lemma}{LEMMA}%[section]
\newtheorem{Theorem}[Lemma]{Theorem}
\newtheorem{Proposition}[Lemma]{Proposition}
\newtheorem{Corollary}[Lemma]{Corollary}

\newtheorem{Example}[Lemma]{Example}
\newtheorem{definition}[Lemma]{Definition}
\newtheorem{exercise}[Lemma]{Exercise}

\newtheorem{remark}[Lemma]{Remark}

\def\bt{\begin{Theorem}}
\def\et{\end{Theorem}}
\def\bp{\begin{Proposition}}
\def\ep{\end{Proposition}}
\def\bcor{\begin{Corollary}}
\def\ecor{\end{Corollary}}
\def\bl{\begin{Lemma}}
\def\el{\end{Lemma}}

\def\beg{\begin{Example}}
\def\eeg{\end{Example}}

\def\bc{\begin{center}}
\def\ec{\end{center}}

\def\vsq{\vskip .25cm}

\def\noi{\noindent}

\def\beq{\begin{equation}}
\def\eeq{\end{equation}}

\def\beqarray{\begin{eqnarray*}}
\def\eeqarray{\end{eqnarray*}}

\def\<{\langle}
\def\>{\rangle}
\def\({\left(}
\def\){\right)}
\def\[{\left[}
\def\]{\right]}

\def\q{\quad}
\def\h{\hbox}

\def\a{\alpha}
\def\b{\beta}
\def\g{\gamma}
\def\d{\delta}
\def\r{\rho}

\def\l{\lambda}

\def\t{\tau}

\def\f{\varphi}

\def\R{\mathbb{R}}

\def\bi{\begin{itemize}}
\def\i{\item}
\def\ei{\end{itemize}}

\def\bpf{\begin{proof}}
\def\epf{\end{proof}}

\def\bq{\begin{quote}}
\def\eq{\end{quote}}

\def\ben{\begin{enumerate}}
\def\een{\end{enumerate}}

\def\bit{\begin{itemize}}
\def\eit{\end{itemize}}

\def\bedef{\begin{definition}}
\def\endef{\hfill $\lhd$\end{definition}}

\def\bex{\begin{exercise}}
\def\eex{\hfill $\lhd$\end{exercise}}
\def\ds{\displaystyle}

\def\brem{\begin{remark}}
\def\erem{\hfill $\lhd$\end{remark}}

\begin{document}

%\hfill M.T. Nair May 20, 2013

\title[On truncated spectral regularization]{On truncated spectral regularization for an \\ ill-posed evolution equation}
\author{M.Thamban  Nair}
\address{Department of Matheamtics, IIT Madras, Chennai, INDIA}
\email{mtnair@iitm.ac.in}
\maketitle

\begin{abstract} In this note we consider the {\it spectral truncation} as the regularization for an ill-posed non-homogeneous parabolic final value problem, and obtain error estimates under a genral source condition when the data, which consist of the non-homogeneous term as well as the final value, are noisy. The resulting error estimate is compared with the corresponding estimate under the Lavrentieve method, and showed that the truncation method has no index of saturation.
\end{abstract}



\section{Introduction} 
Let $H$ be a  Hilbert space and  $A: D(A)\subseteq H\to H$ be a  densely defined positive self adjoint (unbounded) operator.
Given $\f_0\in H$ and  $f\in L^1([0, \infty), H)$,  consider the initial value problem (IVP): 
\beq\label{IVP-1}    \frac{d}{dt}u(t) + Au(t) = f(t),\q u(0)=\f_0.\eeq
It is known \cite{pazy}  that 
if $u(\cdot)$ is a solution for (\ref{IVP-1}), then  
\beq\label{mild-sol-IVP}   u(t) = e^{-tA}\f_0 + \int_0^t e^{-(t-s)A}f(s) ds.\eeq
Here, for $r>0$, the operator $e^{-rA}$ is defined by 
$$e^{-rA}\f := \int_0^\infty e^{r\l} dE_\l \f,\q \f\in H,\q \f\in H.$$
We may also recall (see, e.g., \cite{Yosida}) that, corresponding to any $r\in \R$, the operator 
$e^{rA}$ is defined by 
$$e^{rA}\f := \int_0^\infty e^{r\l} dE_\l \f,\q \f\in H,\q \f\in D(e^{rA}),$$
where 
$$D(e^{rA}):=\{\f\in H: \int_0^\infty e^{2r\l} d\|E_\l\f\|^2<\infty \}.$$
In fact, by spectral theorem, for  any continuous function $g:[0, \infty)\to \R$,
$$g(A)\f:= \int_0^\infty g(\l) dE_\l\f, \q \f\in D(g(A)), $$
where 
$$D(g(A)):=\{\f\in H: \int_0^\infty |g(\l)|^2d\|E_\l\f\|^2 <\infty\},$$
and $g(A)$ is a  positive self adjoint operator. we have  
$$\|g(A)\f\|^2 := \int_0^\infty |g(\l)|^2  d\|E_\l\f\|^2.$$
In particular,  $t\geq 0$
$$e^{tA}\f:= \int_0^\infty e^{t\l} dE_\l\f,\q \f\in D(e^{tA}),$$
where
$$D(e^{tA}):=\{\f\in H: \int_0^\infty e^{2t\l}d\|E_\l\f\|^2<\infty\}$$
and 
$$\|e^{tA}\f\|^2 = \int_0^\infty e^{2t\l}d\|E_\l\f\|^2\} \geq \|\f\|^2.$$
Thus, $e^{tA}$ is a self adjoint operator which is also bounded below, so that it is onto, and hence  
\ben
\i $e^{tA}$ is one-one, onto, and has bounded inverse, namely,  $e^{-tA}$;

\i\, $R(e^{-tA}) = D(e^{tA})\q\forall\, t\geq 0$;

\i ${\mathcal S}:=\{S(t):= e^{-tA}: t\geq 0\}$ is a strongly continuous semigroup on $H$ with  
$$\Big\|\frac{S(t)\f -\f}{t} +A\f\Big\| \to 0\q\h{as}\q t\to 0.$$
\een
Further\,(see \cite{pazy}), $-A$ is the infinitesimal generator of ${\mathcal S}$, i.e.,
$$-A\f:= \lim_{t\to 0}\frac{S(t)\f -\f}{t},\q \f\in D(-A)$$
and $\ds D(-A):=\{\f\in H: \lim_{t\to 0}\frac{S(t)\f -\f}{t} \h{ exists}\}.$

In this paper, we are inerested in the final value problem (FVP), that is the problem of solving 
\begin{eqnarray} \frac{d}{dt}u(t) + Au(t) &=& f(t),\q 0\leq t <\t, \label{FVP-1}\\
u(\t) &=& \f_\t,\label{VP-1}\end{eqnarray} for a known $\f_\t$ for some $\t>0$.
 
Suppose  $u$ is a solution of (\ref{FVP-1})-(\ref{VP-1}). Then, from (\ref{mild-sol-IVP}), we have 
\beqarray
\f_\t &=& e^{-\t A}\f_0 + \int_0^t e^{-(\t-s)A}f(s) ds + \int_t^\t e^{-(\t-s)A}f(s) ds \\
&=& e^{-(\t-t) A}\Big[ e^{-t A}\f_0 + \int_0^t e^{-(t-s)A}f(s) ds + \int_t^\t e^{-(t-s)A}f(s) ds\Big]\\
&=& e^{-(\t-t) A}\Big[u(t) +  \int_t^\t e^{-(t-s)A}f(s) ds\Big]\\
\eeqarray
Thus,
\beq\label{mild-sol-FVP}  u(t) = e^{(\t-t) A}\Big[\f(\t) - \int_t^\t e^{-(\t-s)A}f(s) ds\Big].\eeq
Since the operator $e^{(\t-t) A}$ is unbounded, the above  representation of $u(t)$ shows that 
\bq
{\it small error in the data $(\f_\t, f)$ can lead to large error in the solution $u(t)$}.\eq
In other words, the problem of solving the FVP (\ref{FVP-1})-(\ref{VP-1}) is {\it ill-posed}. 

\bedef\label{def-milde} \rm
If $\f_\t\in H$ and $f\in L^1([0, \t], H)$ are such that
$$\psi(t):= \f_\t - \int_t^\t e^{-(\t-s)A}f(s) ds\q\forall\, t\in [0, \t)$$
belongs to $D(e^{(\t-t) A})$, then $u(\cdot)$ defined by
$$   u(t) = e^{(\t-t)A}\psi(t)$$
is called the {\bf mild solution} of the FVP (\ref{FVP-1}).
\endef

It is to be  observed that a mild solution of the FVP need not be  a solution of the FVP.  In fact, we have the following characterization of the solution of the FVP.

\bt\cite{jana}
Let $\phi\in D(e^{\t A})$ and let $u: [0, \t]\to H$ be defined by $u(t)=e^{(\t-t) A}\phi$, $t\geq 0$. Then $u$  is  a solution of  the FVP
$$u_t + Au(t) = f(t),\q u(\t)=\phi$$
 if and only if $\phi\in D(  Ae^{\t A})$.
\et

\bpf
Note that, for $t\geq o$ and $h>0$, 
$$\frac{u(t+h)-u(t)}{h}= \frac{e^{(\t-t-h) A}\phi-e^{(\t-t) A}\phi}{h} = \frac{e^{-h A}u(t)-u(t)}{h}  .$$
Since $-A$ is the infinitesimal generator of the semigroup $\{e^{-h A}: h \geq 0\}$, it follows that 
$$\lim_{h\to 0}\frac{u(t+h)-u(t)}{h}\h{ exists iff  } u(t)\in D(-A)\h{  iff  } \phi\in D( A e^{\t A}) \q\forall\, t\geq 0,$$
and in that case $u'(t) = -A u(t)$ $u(\t)=\phi$.
\epf

In view of the representation (\ref{mild-sol-FVP}) of the mild solution $u(t)$ of the FVP, the problem of finding $u(t)$  with $u(\t)=\f_\t$ can be posed as a problem of solving the {\it ill-posed  operator equation}
\beq   
\label{op-eq}{\mathcal A}_{t} u(t) =  \psi(t),\eeq
where ${\mathcal A}_{t}: H\to H$ is the bounded operator defined by 
$${\mathcal A}_{t} \f:=  e^{-(\t-t)A}\f,\q \f\in H,$$
and $\psi(t):= \f_\t - \int_t^\t e^{-(\t-s)A}f(s) ds.$
Note that
\bit
\i ${\mathcal A}_{t}$ is an injective bounded self adjoint operator,

\i $R({\mathcal A}_{t}) = D(e^{(\t-t)A})$ is dense in $H$, and 

\i ${\mathcal A}_{t}^{-1}= e^{(\t-t)A}: R({\mathcal A}_{t})\to H$ is not continuous.
\eit

In order to obtain stable approximation for the mild solution given in (\ref{mild-sol-FVP}) for the  FVP, we shall apply the so called {\it  truncated spectral regularization}, and obtain error estimate under a {\it general source condition}. The obtained rate will be compared with the rate resulting from the Lavrentieve method\,\cite{nair-taut}.

\section{Truncated Spectral Regularization (TRS)}

For $\f_\t$ and $f\in L^1([0, \t], H)$, the mild solution, as in (\ref{mild-sol-FVP}), of the FVP has the spectral representation
\beq\label{sol-mild-1}   u(t)  = \int_0^\infty e^{(\t-t)\l}dE_\l(\psi(t))\eeq
whenever 
$$\psi(t):=  \f_\t - \int_t^\t e^{-(\t-s)A}f(s) ds$$ 
belongs to 
$$D(e^{(\t-t) A}):=\{\f\in H:  \int_0^\infty e^{2(\t-t)\l}d\|E_\l\f\|^2<\infty\}.$$ The above representation (\ref{sol-mild-1})  involving integral over the whole of $[0, \infty)$ suggests that a truncation of the inteegral would give a reasonable approximation for $u(t)$. That is the idea in {\it truncated spectral regularization}\,\cite{Tuan, jana, jana-nair}. 

\bedef\label{def-tsr} \rm
The {\it truncated spectral regularized solution}  for the mild solution is defined by
\beq\label{sol-mild-reg-1}    u_\b(t) =  \int_0^\b e^{(\t-t)\l}dE_\l(\psi(t))\eeq
for each $\b>0$. 
\endef

The following theorem shows that $u_\b(\cdot)$ is  an approximation of $u(\cdot)$ for large $\b$.

\bt Under the assumption  $\psi(t)  \in D(e^{\t A})$, 
$$\|u(t)-u_\b(t)\|\to 0\q\h{as}\q \b\to\infty.$$
\et 

\bpf
Since $$\|u(t)\|^2 =  \int_0^\infty  e^{2(\t-t)\l}d\|E_\l(\psi(t))\|^2<\infty,$$
we obtain
$$\|u(t) - u_\b(t)\|^2 =  \int_\b^\infty  e^{2(\t-t)\l}d\|E_\l(\psi(t))\|^2\to 0\q\h{as}\q \b\to\infty.$$
\epf
  
Next, we show that $u_\b(\cdot)$ is, in fact,  stable under perturbations in the data $(\f_\t, f)$.

Suppose  $\tilde \f_\t\in H$ and $\tilde f\in L^1([0, \t], H)$ are the noisy data, in place of the actual data  $\f_\t$ and $f$, respectively.
Let 
$$\tilde u_\b(t) =  \int_0^\b e^{(\t-t)\l}dE_\l(\tilde \psi(t)),$$
where 
$$\tilde \psi(t):=  \tilde \f_\t - \int_t^\t e^{-(\t-s)A}\tilde f(s) ds. $$

\bt\label{Th-reg-1} Let $\f_\t, \tilde \f_\t\in H$ and   $f, \tilde f\in L^1([0, \t], H)$. The for each $t\in [0, \t]$ and  $\b>0$,
\beqarray
\|u_{\b}(t)-\tilde u_{\b}(t)\| &\leq &  e^{(\t-t)\b}\|\psi(t)-\tilde \psi(t)\|\\
&\leq & e^{(\t-t)\b}( \|\f_\t-\tilde\f_\t\|+\|f-\tilde f\|_1 ).
\eeqarray
\et

\bpf 
We observe that 
$$u_\b(t) - \tilde u_\b(t) =  \int_0^\b e^{(\t-t)\l}dE_\l(\psi(t)-\tilde \psi(t))$$
so that
\beqarray
\|u_\b(t) - \tilde u_\b(t)\|^2 
&=&  \int_0^\b e^{2(\t-t)\l}dE_\l\|\psi(t)-\tilde \psi(t)\|^2\\ 
&\leq &  e^{2(\t-t)\b}\int_0^\b dE_\l\|\psi(t)-\tilde \psi(t)\|^2\\
&\leq &  e^{2(\t-t)\b} \|\psi(t)-\tilde \psi(t)\|^2.
\eeqarray
Note that
$$\psi(t)-\tilde \psi(t) = \f_\t-\tilde \f_\t) - \int_t^\t e^{-(\t-s)A}\tilde (f(s)-\tilde f(s)) ds$$
so that  
\beqarray
\|\psi(t)-\tilde \psi(t)\| 
&\leq & \|\f_\t-\tilde \f_\t\| +  \int_t^\t \| e^{-(\t-s)A}\|\,\| f(s)-\tilde f(s)\| ds\\
&\leq & \|\f_\t-\tilde \f_\t\| +  \int_t^\t  \|f(s)-\tilde f(s)\| ds\\
&\leq & \|\f_\t-\tilde \f_\t\| +  \|f-\tilde f\|_1.
\eeqarray
Thus, we obtain the required result.
\epf

We see that the map 
$$(\f, f)\mapsto \|(\f, f)\|_*:= \|\f\|+\|f\|_1$$ 
define a norm on $H\times L^1([0, \t], H)$. Thus, Theorem \ref{Th-reg-1} shows that the truncated spectral regularized solution $u_\b(t)$ is stable under perturbations in the data $(\f_\t, f)$ with respect to the  above norm $\|\cdot\|_*$.

\section{Convergence and error estimates}

\subsection{Convergence} 
From Theorem \ref{Th-reg-1}, the following theorem is immediate.

\bt Let $\f_\t, \tilde \f_\t\in H$ and   $f, \tilde f\in L^1([0, \t], H)$ be such that 
$$\|\f_\t-\tilde\f_\t\|+\|f-\tilde f\|_1\leq \d$$
for some $\d>0$.
Then for each $t\in [0, \t]$ and  $\b>0$,
$$
\|u(t) -  \tilde u_{\b}(t)\| \leq \|u(t)-u_\b(t)\|+ 
e^{(\t-t)\b}\d.$$
Further,  if 
$$\b\approx  \frac{1}{\t-t} \log\Big(\frac{1}{\d^p}\Big)$$
for $0<p<1$, then $e^{(\t-t)\b}\d=\d^p$ and 
$$
\|u(t) -  \tilde u_{\b}(t)\| = o(1)\q\h{as}\q \d\to 0.$$
\et

\subsection{Estimates under general source condition}

% 
% $$u(t)-  u_\b(t) =  \int_\b^\infty  e^{(\t-t)\l}dE_\l(\psi(t)) $$
%
%Suppose $u(t) = e^{-(\t-t)A}v(t)$. Then 
%$$v(t) = e^{(\t-t)A} u(t) = e^{2(\t-t)A} \psi(t) $$
%\beqarray
%\|u(t)-  u_\b(t)\|^2 &=&   \int_\b^\infty  e^{2(\t-t)\l}dE_\l\|\psi(t)\|^2 \\
%& =&   \int_\b^\infty  e^{-2(\t-t)\l}e^{4(\t-t)\l}dE_\l\|\psi(t)\|^2 \\
%& \leq &   e^{-2(\t-t)\b}\int_\b^\infty  e^{4(\t-t)\l}dE_\l\|\psi(t)\|^2 \\
%& \leq  &   e^{-2(\t-t)\b}\r_\b(t)^2,
%\eeqarray
%where $\r_\b(t)\to 0$ and $\b\to 0$ and $\r_\b(t)\leq \|v(t)\|$.
% 
For obtaining  error  estimates it is required to assume certain smoothness assumptions on the solution, the so called {\it source conditions}. For this purpose, we consider a general condition of the form 
\beq\label{source-1} u(t)\in D(h_t(A)),\eeq
where the function $h_t: [0, \infty)\to (0, \infty)$ is continuous and 
for each $t\in [0, \t)$,
$$h_t(\l)\to \infty\q\h{as}\q \l\to\infty.$$
Note that the condition (\ref{source-1}) is equivalent to 
\beq\label{source-2} \int_0^\infty [h_t(\l)]^2d\|E_\l u(t)\|^2<\infty.\eeq
At this point one may recall that in \cite{Tuan}, Tuan has considered the source conditions of the forms,
\beq\label{source-3}  \int_0^\infty \l^{2p} d\|E_\l u(t)\|^2<\infty \q\h{and}\q \int_0^\infty e^{2\l q} d\|E_\l u(t)\|^2<\infty,\eeq
for $p, q>0$ 
%yielding estimates 
%$$\|u(t)-u_\b(t)\| = O(e^{-t\b}\b^{-p})\q\h{and}\q 
%\|u(t)-u_\b(t)\| = O(e^{-t(\b+q)}),$$
respectively.
Note that the source conditions in (\ref{source-3}) are special cases of (\ref{source-2}) obtained by the  choices 
$$h_t(\l):= \l^p\q\h{and}\q h_t(\l):= e^{q\l},$$
respectively.

\bt\label{Th-reg-2}
Suppose $\f_\t\in H$ and $f\in L^1([0, \t], H)$ are  such that 
$\psi(t)\in D(e^{\t A})$ for each $t\in [0, \t)$ 
and  there exists a monotonically increasing continuous function $h_t: [0, \t]\to (0, \infty)$ such that
\vsq

\ben
\i[\rm(i)] $h_t(\l)\to \infty$ as $\l\to\infty$,
\vsq
\i[\rm(ii)] $u(t)\in D(h_t(A))$, 
\een
Let $\r_t>0$ be such that  $\|h_t(A)u(t)\|\leq \r_t$.
Then
$$\|u(t) - u_\b(t)\|
\leq   \frac{\r_{t,\b}}{h_t(\b)},$$
where $\r_{t,\b}\leq \r_t$ and $\r_{t,\b}\to 0$ as $\b\to\infty$. 
\et

\bpf  Recall from Definition \ref{def-milde} that $u(t) = e^{(\t-t)A}\psi(t).$
Hence, 
\beqarray
\|u(t) - u_\b(t)\|^2 &=&  \int_\b^\infty  e^{2(\t-t)\l}d\|E_\l(\psi(t))\|^2\\
&=&  \int_\b^\infty \frac{1}{[h_t(\l)]^2} [h_t(\l)]^2 e^{2(\t-t)\l}d\|E_\l(\psi(t))\|^2.
\eeqarray
Since $h_t$ is monotonically increasing, from the above, we obtain 
\beq\label{eqn-12} \|u(t) - u_\b(t)\|^2\leq   \frac{1}{[h_t(\b)]^2}\int_\b^\infty  [h_t(\l)]^2 e^{2(\t-t)\l}d\|E_\l(\psi(t))\|^2.\eeq
By the assumption,
\beqarray
\int_0^\infty [h_t(\l)]^2e^{2(\t-t)\l} d\|E_\l(\psi(t))\|^2 
&=& \|h_t(A)e^{(\t-t)A}\psi(t)\|^2 \\
& = &  \|h_t(A)u(t)\|^2\leq \r_t^2.
\eeqarray
Hence, taking $\r_{t, \b}$ such that 
$$\r_{t,\b}^2 =\int_\b^\infty  h_t(\l)^2 e^{2(\t-t)\l}d\|E_\l(\psi(t))\|^2,$$
the inequality (\ref{eqn-12}) leads to   
$$\|u(t) - u_\b(t)\|
\leq   \frac{\r_{t,\b}}{h_t(\b)},$$
where $\r_{t,\b}\leq \r_t$ and $\r_{t,\b}\to 0$ as $\b\to\infty$. 
\epf

\brem\rm 
Recently, Jana \cite{jana} and Jana and Nair \cite{jana-nair}  used similar general source condition, but based on  the data $\f_\t, f$ instead of the mild solution $u(t)$. 
\erem

Combining the last two theorems,
we obtain the following.

\bt\label{Th-est-main-1} 
Suppose $\tilde \f_\t$ and $\tilde f$ are noisy data such that
$$\|\f_\t-\tilde\f_\t\|+\|f-\tilde f\|_1\leq \d$$
for some noise level $\d>0$. Then
$$\|u_{\b}(t)-\tilde u_{\b}(t)\| \leq  e^{(\t-t)\b}\d.$$
If  $\r_{t, \b}$ and $h_t(\cdot)$ are as in Theorem \ref{Th-reg-2},
then we have
$$\|u(t) - \tilde u_\b(t)\|
\leq   \frac{\r_{t,\b}}{h_t(\b)} + e^{(\t-t)\b}\d.$$
\et

\section{Parameter Choice Strategy}

In Theorem \ref{Th-est-main-1}, we obtained an  estimate for the error  
$\|u(t) - \tilde u_\b(t)\|$ a smoothness assumption on $u(\cdot)$.
Now, we choose $\b:=\b_t^\d$ depending on $t$ and $\d$ such that 
the obtained estimate converges to $0$ as $\d\to 0$.

\bt
or $\l>0$, let 
$\xi_t(\l):= h_t(\l)e^{(\t-t)\l}$, and let 
$$\b=\b_t(\d):= \xi_t^{-1}(\r/\d).$$
Then
$$\|u(t) - \tilde u_\b(t)\|
\leq    \frac{2\r}{h(\xi_t^{-1}(\r/\d))} .$$
In particular,
$$\|u(t) - \tilde u_\b(t)\| \to 0\q\h{as}\q \d \to 0.$$
\et

\bpf
Note that
\beqarray
\frac{\r_t}{h_t(\b)}  = e^{(\t-t)\b}\d 
&\iff& \xi_t(\b):= h_t(\b)e^{(\t-t)\b} = \frac{\r_t}{\d} \\
&\iff & \b = \xi_t^{-1}(\r_t/\d).
\eeqarray
Thus, for the choice of $\b = \xi_t^{-1}(\r_t/\d)$, 
Theorem \ref{Th-est-main-1} implies 
\beqarray
\|u(t) - \tilde u_\b(t)\|
&\leq & \frac{\r_t}{h_t(\b)} + e^{(\t-t)\b}\d\\
&\leq&    \frac{2\r}{h(\xi_t^{-1}(\r/\d))}\\
\eeqarray
Since $h(\xi_t^{-1}(\r_t/\d))\to \infty$ {as} $\d\to 0$.
 $$\|u(t) - \tilde u_\b(t)\|\to 0 \q\h{as}\q \d\to 0.$$
\end{proof}

\subsection{A special choice of the source condition} 
\vsq
Suppose 
$u(t) \in R(e^{-\g(\t-t)A})$ for some $\g>0$. 
Then 
$$u(t) = e^{-\g(\t-t)A} v_\g(t)$$ 
for some $v_\g(t)\in H$. Thus, 
$$u(t) \in D(e^{\g(\t-t)A})$$ and 
$$\|e^{\g(\t-t)A}u(t)\|\leq \r_t.$$
Note that the funtion 
$$h_t(\l):= e^{\g(\t-t)\l},\q\l\geq 0,$$
satisfies the properties (1)-(3) in Theorem \ref{Th-reg-2}.
Thus, by Theorem \ref{Th-est-main-1}, 
$$\|u(t) - \tilde u_\b(t)\|
\leq   \r_{t,\b} e^{-\g(\t-t)\b} + e^{(\t-t)\b}\d.$$

\bt
Suppose $u(t) \in R(e^{-\g(\t-t)A})$ for some $\g>0$ and 
$$\|e^{\g(\t-t)A}u(t)\|\leq \r_t. $$
Then 
\beq\label{est-special-1} \|u(t) - \tilde u_\b(t)\|
\leq   \r_{t} e^{-\g(\t-t)\b} + e^{(\t-t)\b}\d.\eeq
Further, taking  
$$\b_t :=  \frac{1}{(\g+1)(\t-t)} \log\Big(\frac{1}{\d}\Big),$$
we have
$$\|u(t) - \tilde u_{\b_t}(t)\|
\leq   (1+\r_{t})\Big[ \log\Big(\frac{1}{\d}\Big) \Big]^{-\g}.$$
\et

\section{Comparison with Lavrentiev regularization}

Recall that  
the operator  ${\mathcal A}_{t}: H\to H$  defined by 
$${\mathcal A}_{t} \f:=  e^{-(\t-t)A}\f,\q \f\in H$$
 is   injective,   continuous,  self adjoint,   with   $R({\mathcal A}_{t})$  dense in $H$.
Let $u(t)$ be the solution of (\ref{op-eq-1}), that is, 
$${\mathcal A}_{t} u(t) = \psi(t):= \f_\t - \int_t^\t e^{-(\t-s)A}f(s) ds.$$
Let  $u^L_\a(\cdot)$ be the Lavrentive regularized solution, i.e., 
$$   ({\mathcal A}_{t} +\a I)u^L_\a(t) = \psi(t) .$$
Then, from the standard theory \cite{ehn, nair-taut},  we know that
 $$\|u(t)- u^L_\a(t)\|\to 0\q\h{as}\q \a\to 0$$
 and 
 \beq\label{est-Lav-1} \|u^L_\a(t)- \tilde u^L_\a(t)\|\leq  \frac{\d}{\a} .\eeq 
Next, suppose 
 \beq\label{source-0} u(t) = {\mathcal A}_{t}^\g v(t)\q\h{with}\q \|v(t)\|\leq \r_t,\eeq
for some $\g>0$.  Equivalently, 
 \beq\label{source-1} u(t)\in D(e^{\g(\t-t)A})\q\h{with}\q \|e^{\g(\t-t)A}u(t)\|\leq \r_t.\eeq

Then we have the estimate 
\beq
\label{est-Lav-1} \|u(t)- u^L_\a(t)\|\leq \r_t \a^\g\q\h{whenever}\q 0<\g\leq 1,\eeq
and consequently,
$$ \|u(t)- \tilde u^L_\a(t)\|\leq \r_t \a^\g +\frac{\d}{\a}.$$
Note that 
\beq\label{parameter-lav} \frac{\d}{\a} = e^{(\t-t)\b} \iff 
\b:= \frac{1}{\t-t}\ln\Big(\frac{1}{\a} \Big).\eeq 
Thus, the estimate in (\ref{est-Lav-1}) is same as (\ref{est-special-1}) for the choice of $\b$ as in (\ref{parameter-lav}).
However, the estimate in (\ref{est-special-1}) is valid for all $\g>0$, whereas (\ref{est-Lav-1}) is valid only for $0<\g\leq 1$.

\vsq
\noi
{\bf Acknowledgements:}\, This work is completed while the author was a visiting mathematician at Sun Yat-sen University, Guanzhou, China, during the period June 13 to July 8, 2019. The support and the warm hospitality received from Prof. Hongqi Yang  are gratefully acknowledged.

\end{document}